\documentclass[12pt]{article}
\usepackage{amsthm,array,amsmath,graphics,amssymb}

\newcommand\Z{{\mathbb Z}}
\newcolumntype{L}{>{$}l<{$}}
\newcolumntype{C}{>{$}c<{$}}

\theoremstyle{plain}
\newtheorem{thm}{Theorem}
\newtheorem{defn}{Definition}

\title{Bandwidth reduction in rectangular grids}

\author{Titu Andreescu, Walter Stromquist and Zoran \v Suni\'k\\
\small American Mathematics Competitions\\[-0.8ex]
\small 1740 Vine St. \\[-0.8ex]
\small Lincoln, NE 68588-0658, USA\\[-0.8ex]
\small \texttt{titu@amc.unl.edu}\\
\small 132 Bodine Road\\[-0.8ex]
\small Berwyn, PA 19312\\[-0.8ex]
\small \texttt{walters@chesco.com}\\
\small Department of Mathematics and Statistics\\[-0.8ex]
\small 810 Oldfather Hall, University of Nebraska\\[-0.8ex]
\small Lincoln, NE 68588-0323, USA\\[-0.8ex]
\small \texttt{zsunik@math.unl.edu}}

\date{
\small MR Subject Classifications: 05C78}

\begin{document}
\maketitle

\begin{abstract}
We show that the bandwidth of a square two-dimensional grid of
arbitrary size can be reduced if two (but not less than two)
edges are deleted. The two deleted edges may not be chosen
arbitrarily, but they may be chosen to share a common endpoint
or to be non-adjacent.

We also show that the bandwidth of the rectangular $n \times m$
($m \geq n$) grid can be reduced by $k$, for all $k$ that are
sufficiently small, if $m-n+2k$ edges are deleted.
\end{abstract}

\section{Introduction} We consider only simple undirected graphs
(no loops, no multiple edges). A \emph{vertex numbering} of a
graph $G=(V,E)$ is a bijective map $\nu:V \to [k]$ from the
vertex set $V$ of $G$ to the set of the first $k$ positive
integers $[k]=\{1,2,\dots,k\}$, where $k=|V|$. The absolute
value of the difference between the numbers at the two endpoints
of an edge $e$ in $E$ is denoted by $f_\nu(e)$ and called the
\emph{length} of $e$ induced by $\nu$. Thus any vertex numbering
of $G$ induces an edge labelling $f_\nu:E \to \Z_+$ of $G$ by
positive integers. The largest length of an edge in $E$ (i.e.
the largest label used by $f_\nu$) is called the \emph{bandwidth
of the numbering} $\nu$. If the edge set is empty the bandwidth
is 0 by definition. The smallest possible bandwidth, taken over
all possible vertex numberings of $G$, is called the
\emph{bandwidth} of the graph $G$. Sometimes the adjective
``linear'' is used due to the following physical interpretation.
Every vertex numbering $\nu$ corresponds to a \emph{linear
arrangement} of the graph $G$ in which the vertex $v$ in $V$ is
placed at the number $\nu(v)$ on the real line. The (linear)
bandwidth of a linear arrangement is then just the length of the
longest wire needed to assemble the graph $G$ and the minimal
such bandwidth, taken over all possible linear arrangements of
$G$, is the (linear) bandwidth of $G$.

The two-dimensional \emph{rectangular grid} $G_{m,n}$, for $m
\geq n \geq 1$, is the graph whose vertices are the points in
the set $V=[m]\times[n]=\{(i,j)|1 \leq i \leq m, 1 \leq j \leq
n\}$ with an edge between two vertices if and only if the
Euclidean distance between them is 1 (think of the natural
embedding of the vertex set in the Cartesian plane). We write
$G_n$ for the square grid $G_{n,n}$. It is well known that the
bandwidth of the rectangular grid $G_{m,n}$ is $n$, unless
$m=n=1$.

\begin{thm}[J. Chv\'atalov\'a
\cite{chvatalova:bandgrid}]\label{thm:chvatalova} Let $m \geq n
\geq 1$. If $m \geq 2$, then the bandwidth of the
two-dimensional rectangular grid $G_{m,n}$ is $n$.
\end{thm}

The bandwidth is, in some sense, a measure of how difficult it
is to embed the graph on a line.  The following result provides
some refinement of this measure for the rectangular grids.

\begin{thm}[P. Fishburn and P. Wright
\cite{fishburn-w:2(n-1)}]\label{thm:2(n-1)} Let $m \geq n \geq
1$. Every vertex numbering of $G_{m,n}$ of bandwidth $n$ induces
at least
\[ 2(n-1) + n(m-n) \]
edges of length $n$. Moreover, the \emph{down diagonal
lexicographic linear arrangement} induces exactly
$2(n-1)+n(m-n)$ edges of length $n$.
\end{thm}

As an example, the down diagonal lexicographical linear
arrangement on the square grid $G_n$ is given by
\[
\begin{tabular}{LLLLLL}
    t+1    &        & \ddots &         & n^2   \\
    \ddots & t+2    &        & \ddots  &       \\
    4      & \ddots & \ddots &         & \ddots\\
    2      & 5      & \ddots & \ddots  &       \\
    1      & 3      & 6      & \ddots  & t+n
\end{tabular},
\]
where $t=(n-1)n/2$. This numbering has bandwidth $n$ and induces
exactly $2(n-1)$ edges of length $n$ (all horizontal edges
incident to the main diagonal vertices).

Note that the down diagonal lexicographical vertex numbering of
$G_n$ that we just displayed indicates that we sometimes think
of $G_{m,n}$ as a rectangular board of dimension $n \times m$
(meaning $n$ rows and $m$ columns) in which the squares
represent the vertices and any pair of squares that have common
side are considered to be neighbors.

By the result of Chv\'atalov\'a,  we know that in order to make
a linear arrangement of the square grid $G_n$ we must be ready
to use pieces of wire of length at least $n$. Moreover, by the
result of Fishburn and Wright, in order to make an arrangement
that does not use any pieces longer than $n$, we must be ready
to use at least $2(n-1)$ pieces of length $n$.

But what if we do not have such long pieces? In such a case the
most practical thing one can do is to try to assemble as large
part of the graph as possible, which amounts to deletion of some
of the edges.

\section{Bandwidth reduction} The bandwidth of the path $P_n$ of
length $n-1 \geq 1$ is 1, of the cycle $C_n$ of length $n \geq
3$ is 2, and of the complete graph $K_n$ on $n \geq 1$ vertices
is $n-1$. The deletion of any edge in the cycle $C_n$ produces
the path $P_n$ and thus reduces the bandwidth by 1. Similarly,
the deletion of any edge in $K_n$ also reduces the bandwidth by
1. We are interested in the minimal number of edges that need to
be deleted in an arbitrary graph in order to reduce its
bandwidth.

\begin{defn}
The \emph{bandwidth reduction number} of a graph $G$ of
bandwidth $b$, $b \geq 1$, is the minimal number of edges that
need to be deleted from $G$ in order to obtain a subgraph of
bandwidth no greater than $b-1$.
\end{defn}

For example, cycles and complete graphs have bandwidth reduction
number 1, while paths have bandwidth reduction number equal to
their length. As a more interesting example, consider the wheel
$W_m$, $m \geq 4$, a graph on $m$ vertices that consists of a
cycle of length $m-1$ together with an ``center'' vertex
connected to every vertex in the cycle by an edge. The following
table provides the bandwidth and bandwidth reduction number of
the wheel graphs. We do not provide the easy proofs, but we note
the exceptional case of $W_7$, which is the only wheel whose
bandwidth reduction number is 3.
\[
\begin{tabular}{lCCCCCC}
 graph     & W_4 & W_5 & W_6 & W_7 & W_{2n}\;(n\geq 4) & W_{2n+1}\; (n\geq 4)\\
 \hline
 bandwidth & 3   & 3   &  3  & 3   & n                 & n\\
 bandwidth reduction n.
           & 1   & 1   &  2  & 3   & 1                 & 2
\end{tabular}
\]
Another relatively easy example is provided by the complete
bipartite graphs $B_{m,n}$, for $m \geq n \geq 1$.
\[
\begin{tabular}{lCCC}
 graph     & B_{2,2} & B_{2k+1,n} & B_{2k,n}\;(k \neq 1\text{ or }n \neq 2) \\
 \hline
 bandwidth & 2   & k+n   &  k+n-1\\
 bandwidth reduction number
           & 1   & 1   &  2
\end{tabular}
\]

The notion of bandwidth reduction can be further extended as
follows.

\begin{defn}
Let $G$ be a graph of bandwidth $b$. For $k=1,\dots,b$, the
$k$-th bandwidth reduction number of $G$, denoted by $br_k(G)$,
is the minimal number of edges that need to be deleted from $G$
in order to obtain a subgraph of bandwidth no greater than
$b-k$.
\end{defn}

We note that the $k$-th bandwidth reduction number of a graph
$G$ of bandwidth $b$ is the minimal number of edges $e$ of
induced length $f_\nu(e)$ greater than $b-k$, taken over all
possible numberings $\nu$ of $G$, i.e.,
\[ br_k(G) = \min_{\nu} \left|\{e \in E | f_\nu(e)>b-k\}\right|. \]
For example, in the case of the complete graph $K_n$ we have
\[ br_k(K_n) = 1+2+ \cdots k = \frac{k(k+1)}{2}, \]
for $k=1,\dots,n-1$, since the bandwidth of $K_n$ is $n-1$ and
every numbering of $K_n$ induces $i$ edges of length $n-i$, for
$i=1,\dots,n-1$. This implies that if a graph on $n$ vertices
has more than
\[ \frac{n(n-1)}{2} - \frac{k(k+1)}{2} \]
edges, then its bandwidth is at least $n-k$, for
$k=1,\dots,n-1$.

Before we move on, we observe that the bandwidth of a graph can
be reduced by more than 1 by deletion of a single edge. This is
why we require the newly obtained graph to have bandwidth no
greater than $b-k$ rather than exactly $b-k$ in the definition
of the $k$-th bandwidth reduction number.

Consider the graph $G$ that consists of two copies of the wheel
$W_7$ with an additional edge between the wheel centers (two
wheels with an axis between them). Since an arbitrary graph on
$v$ vertices with diameter $d$ and bandwidth $b$ satisfies
\[ 1+ db \geq v, \]
the bandwidth of the graph $G$ is at least 5 ($v$=14 and $d=3$).
However the deletion of the edge between the centers of the two
wheels leads to a graph of bandwidth 3. Thus $br_1(G)=br_2(G)=1$
for this graph.

We give now an upper bound on the $k$-th bandwidth reduction
number in rectangular grids, for small $k$.

\begin{thm}\label{reduction}
Let $m \geq n > 2k$. The $k$-th bandwidth reduction number of
the rectangular grid $G_{m,n}$ satisfies
\[ br_k(G_{m,n}) \leq m-n+2k. \]
\end{thm}
\begin{proof}
We construct an example of a numbering $\nu$ of $G_{m,n}$ in
which only $m-n+2k$ edges are longer than $n-k$. The example is
a modification of the down diagonal lexicographical linear
arrangement.

Think of $G_{m,n}$ as of rectangular board of dimension $n
\times m$. First modify the board by cutting out the lower right
part of the board of dimension $k \times (m-n+k+1)$ and flipping
it over as indicated in Figure~\ref{fig:reduction}.
\begin{figure}[!ht]
  \begin{center}
  \includegraphics{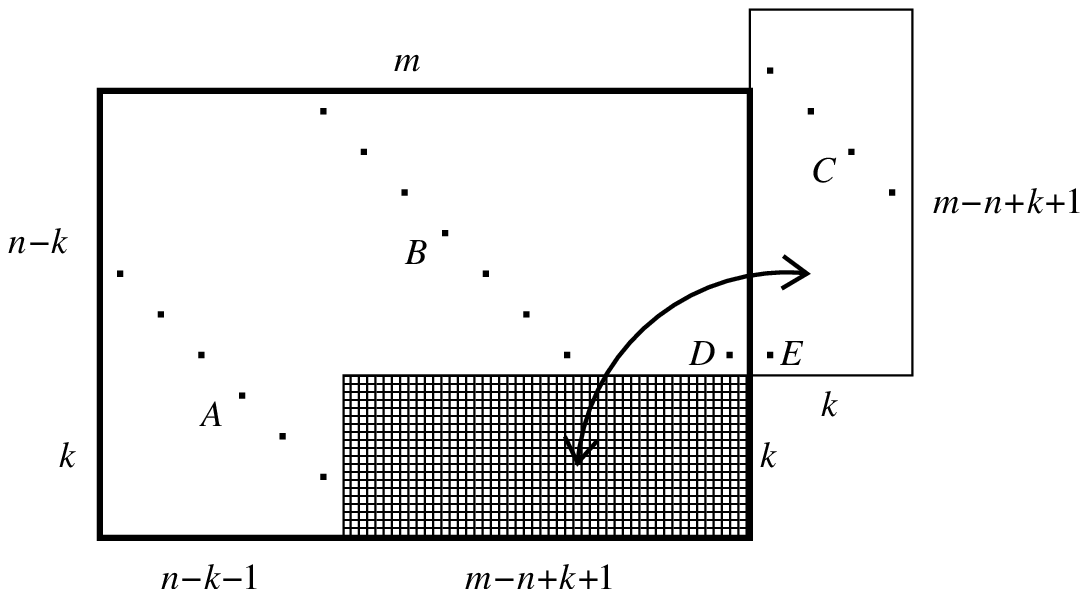}
  \end{center}
  \caption{A numbering of $G_{m,n}$ with $m-n+2k$ edges longer than $n-k$}
  \label{fig:reduction}
\end{figure}

Enumerate all the squares in the modified board in the down
diagonal lexicographical fashion. Consider a ``typical'' vertex
$A$ that lives on a diagonal of length $n-k-1$. The vertical
edges incident to $A$ have length $n-k-1$ and the horizontal
ones have length $(n-k-1)+1=n-k$. A ``typical'' vertex $B$ that
lives on a diagonal of length $n-k$ is incident to horizontal
edges of length $n-k$ and vertical edges of length $(n-k)-1$. A
``typical'' vertex $C$ that lives on a diagonal of length $k$ is
incident to vertical edges of length $k$ and horizontal edges of
length $k+1=2k+1-k \leq n-k$, where the inequality comes from
the assumption that $2k<n$. In all ``atypical'' cases, near the
bottom left or the upper right corner(s) of the modified board,
the diagonals are even shorter, which implies that the incident
edges are also shorter.

Thus the bandwidth of the given numbering of the modified board
is $n-k$. After we flip back the modified part of the board to
its original position the newly created $k$ horizontal edges
(between the non-shaded and the shaded part of the board) will
have length longer than $n-k$ and all but one of the newly
created $m-n+k+1$ vertical edges will have length longer than
$n-k$. Indeed the induced length of the rightmost vertical edge
between the shaded and the non-shaded part is exactly $n-k$,
since its endpoints (denoted by $D$ and $E$ in
Figure~\ref{fig:reduction}) were already neighbors in the
modified board.
\end{proof}

We provide one example to illustrate the preceding result.
Assume that we want to reduce the bandwidth of $G_8$ to 6. Thus
$m=n=8$, $k=2$ and the numbering of the modified board is given
by
\[
\begin{tabular}{LLLLLLLLLLL}
 26 & 32 & 38 & 44 & 50 & 56 & 61 & 64 \\
 21 & 27 & 33 & 39 & 45 & 51 & 57 & 62 \\
 16 & 22 & 28 & 34 & 40 & 46 & 52 & 58 \\
 11 & 17 & 23 & 29 & 35 & 41 & 47 & 53 & 59 & 63\\
  7 & 12 & 18 & 24 & 30 & 36 & 42 & 48 & 54 & 60\\
  4 &  8 & 13 & 19 & 25 & 31 & 37 & 43 & 49 & 55\\
  2 &  5 &  9 & 14 & 20\\
  1 &  3 &  6 & 10 & 15
\end{tabular}
\]

Another way to state Theorem~\ref{reduction} is to say that the
incidence matrix of $G_{m,n}$ can be written in such a way that
all non-zero entries, except for $m-n+2k$ symmetric pairs, are
within distance $n-k$ from the main diagonal.

For example, this means that the incidence matrix of the square
grid $G_n$ can be written in such a way that all non-zero
entries, except for 2 symmetric pairs, are within distance $n-1$
from the main diagonal. This contrasts nicely with
Theorem~\ref{thm:2(n-1)}, which implies that if all non-zero
entries of the incidence matrix of $G_n$ are within distance $n$
from the main diagonal then there are at least $2(n-1)$
symmetric pairs of non-zero entries at distance $n$ from the
main diagonal.

We show now that the bandwidth of square grids cannot be reduced
if only 1 edge is deleted. The exceptional cases are $n=1$ and
$n=2$. Indeed, $G_1$ has bandwidth 0 which cannot be reduced,
while $G_2$ is the cycle $C_4$ and its bandwidth reduction
number is 1.

\begin{thm}\label{thm:2}
Let $n \geq 3$. The bandwidth reduction number of $G_n$ is 2,
i.e., the minimal number of edges that needs to be deleted from
$G_n$ in order to obtain a graph of bandwidth $n-1$ is 2.
Moreover, the edges that need to be deleted may, but do not have
to, share a common endpoint.
\end{thm}

\begin{proof}[Proof of Theorem~\ref{thm:2}]
For an arbitrary numbering $\nu$ of $G_n$, call an edge
\emph{long} if its induced length is at least $n$ and call it
\emph{short} otherwise. An example of a numbering of $G_n$, for
$n \geq 3$, with only two long edges is already provided in the
proof of Theorem~\ref{reduction}. We note that the two long
edges in that example share a common vertex, namely the vertex
in row $n$ column $n-1$. The lengths of the long edges are
$5n-7$ and $3n-4$.

Rather than providing a general example of a numbering with only
two long edges that are not adjacent, we provide an example for
the case $n=6$, from which one can easily extract the general
pattern. The two long edges, which are the two leftmost vertical
edges between the top two rows, have lengths $5n-8$ and $3n-5$.
\[
\begin{tabular}{LLLLLL}
 33 & 29 & 25 & 30 & 34 & 36\\
 11 & 16 & 21 & 26 & 31 & 35\\
 7  & 12 & 17 & 22 & 27 & 32\\
 4  & 8  & 13 & 18 & 23 & 28\\
 2  & 5  &  9 & 14 & 19 & 24\\
 1  & 3  & 6  & 10 & 15 & 20
\end{tabular}.
\]

It remains to show that every numbering of $G_n$, $n \geq 3$,
has at least two long edges.

Let $\nu$ be an arbitrary labelling of $G_n$. Choose the
smallest $k$ such that all rows or all columns of $G_n$ have a
label in $[k]$. Without loss of generality we may assume that
all rows have a label in $[k]$. Define the partial numbering
$\kappa$ of $G_n$ to be the numbering of the vertices obtained
from the numbering $\nu$ by deletion of the labels greater than
$k$ (thus only $k$ vertices have a label).

\textbf{Claim 1.} The label $k$ is alone in its row in the
partial numbering $\kappa$.

Otherwise each row would have a label from $[k-1]$ which
contradicts the minimality in the choice of $k$.

\textbf{Claim 2.} No row is completely numbered by $\kappa$.

Otherwise each column would have a label from $[k-1]$ which
contradicts the minimality in the choice of $k$.

Call an edge between vertex numbered by an element in $[k]$ and
an element not in $[k]$ a \emph{boundary edge}. If there are at
least $n+1$ boundary edges then at least two of them would have
to be long. A direct corollary of Claim~2 is that each row has
at least one horizontal boundary edge. This already means that
there must be at least one long edge (apropos, this proves
Theorem~\ref{thm:chvatalova} in the case of square grids). Thus,
we assume in the sequel that there are exactly $n$ horizontal
boundary edges and only one of them is long. Moreover, we assume
that this is the only long edge in $\nu$ and we seek a
contradiction.

\textbf{Claim 3.} The endpoints of the $n-1$ short horizontal
boundary edges are numbered by the pairs
\[ (k,k+n-1), (k-1,k+n-2), (k-2,k+n-3), \dots, (k-n+2,k+1),  \]
while the endpoints of the long edge come from a pair
$(s,\ell)$, where $s$ is a ``small'' number less than or equal
to  $k-n+1$ and $\ell$ is a ``large'' number greater or equal to
$\geq k+n$.

Indeed, any other numbering of the endpoints of the boundary
edges would result in at least two long horizontal edges.

\textbf{Claim 4.} For every row $i$ in $G_n$, the vertices
labelled by the partial numbering $\kappa$ form a horizontal
path $p_i$ that ends at the leftmost or the rightmost column.

This follows from the fact that every row has exactly one
horizontal boundary edge.

Thus each row contains a unique maximal path that consists of
vertices numbered by $\kappa$. These $n$ paths will be called
$\kappa$-paths.

\textbf{Claim 5.} The $\kappa$-path in any row next to $k$
consist of a single vertex in the column of $k$.

Indeed, if the neighboring row has a vertex labelled by $\kappa$
that is not in the same column as $k$, this vertex would be
involved in a vertical edge that would have to be long (all
vertices in the row of $k$ are labelled by numbers greater or
equal to $k+n-1$).

Without loss of generality we may assume that the vertex $k$
(together with its vertical neighbors from $[k]$) is in the
leftmost column in $G_n$.

\textbf{Claim 6.} If the $\kappa$-path $p$ is no longer than
$n-2$, then the $\kappa$-path $q$ in any neighboring row
occupies the same columns as $p$ or differs in exactly one
column.

Again, in any other case there would be a long vertical edge
involving a number from $[k]$ and a number outside of $[k]$.
Note that we required that the length of $p$ be no longer than
$n-2$, since if the length is $n-1$ the neighboring
$\kappa$-path $q$ may also have length $n-1$ and start at the
opposite end without creating additional long vertical edge(s).
However, we will see that this does not happen.

\textbf{Claim 7.} No horizontal $\kappa$-path is longer than
$n-2$ and they all start at the left end.

The length of the $\kappa$-path $k$ is 1 and so is the length of
the neighboring $\kappa$-path(s). Going further away, the length
of the $\kappa$-paths can only go up by 1 from one row to
another, so the only way we can have a $\kappa$-path of length
$n-1$ is if $k$ is in the bottom or in the top row, the next
$\kappa$-path has length 1 and the length of the remaining
$\kappa$-paths grows by 1 as we move away from the row of $k$.
Without loss of generality assume that $k$ is in the bottom left
corner. Consider now the $\kappa$-path $p_s$, corresponding to
the row with a long horizontal edge, together with the
horizontal path $q$ in the row below.
\[
\begin{tabular}{LLLLL}
* & \dots & * & s & \ell \\
* & \dots & * & x &
\end{tabular},
\]
where the stars represent arbitrary numbers in $[k]$. Clearly,
if $x$ is not in $[k]$ we obtain a long vertical edge with
endpoints $x$ and $s$. Thus $x$ is in $[k]$ and the lengths of
$p_s$ and $q$ are equal. Since the lengths of the $\kappa$-paths
start at 1 and repeat at the row of $k$ and at the row of $s$,
the longest $\kappa$-path can only reach the length $n-2$.

\textbf{Claim 8.} The $\kappa$-path $p_s$ corresponding to the
row with the long horizontal edge $(s, \ell)$ has the same
length as the $\kappa$-paths in the rows next to it.

We already proved that the $\kappa$-path below $p_s$ cannot be
shorter than $p_s$. It cannot be longer as well, since the
configuration
\[
\begin{tabular}{LLLLL}
* & \dots & * & s & \ell \\
* & \dots & * & * & x
\end{tabular},
\]
with $x$ in $[k]$, indicates that $\ell$ and $x$ are the
endpoints of a long vertical edge. By symmetry, the
$\kappa$-path above $p_s$ also has length equal to the length of
$p_s$. This finishes the proof of Claim 8.

We now consider the partial numbering $\tau$ that corresponds to
the set of numbers $[k+n-1]$. As observed before, the numbers
$k+1,k+2,\dots,k+n-1$ are placed in different rows, one in each
row except for the row of the long horizontal edge $s\ell$. Thus
$n-1$ horizontal boundary edges are now formed between the
vertices numbered by $[k+n-1]$ and those not numbered by
$[k+n-1]$. In addition, there exist at least one vertical
boundary edge with endpoints in an element of $[k+n-1]$ and
$\ell$. These $n$ boundary edges of the partial numbering $\tau$
must contain a long edge, a contradiction.
\end{proof}

\section{Final remarks}
We conjecture that in order to reduce the bandwidth of $G_{m,n}$
by $k$, for $k<n/2$, we need to delete at least $m-n+2k$ edges,
i.e., we conjecture that the upper bound provided by
Theorem~\ref{reduction} is the actual value of the bandwidth
reduction number for all sufficiently small values of $k$.
Theorem~\ref{thm:2} shows that this conjecture is true when
$m=n$ and $k=1$.

The proof of Theorem~\ref{thm:2} indicates an easy way to see
that at least $k$ edges need to be deleted in order to reduce
the bandwidth of $G_{m,n}$ by $k$. Namely, the first time some
partial numbering has a representative in each row (or each
column) there would be at least $n$ horizontal (or vertical)
boundary edges at least $k$ of which must have length at least
$n-k+1$.

For a numbering $\nu$ of a graph $G=(V,E)$ one can define the
\emph{vertex-isoperimetric} number of $\nu$ to be the maximal
number of vertices in $G$ numbered by $[k]$ that have neighbors
that are not numbered by $[k]$, where the maximum is taken over
all $k=0,\dots,|V|$. The smallest possible vertex-isoperimetric
number, taken over all numberings of $G$, is the
\emph{vertex-isoperimetric number} of the graph $G$. If $b(G)$
is the bandwidth of $G$ and $vi(G)$ is its vertex-isoperimetric
number then
\[ b(G) \geq vi(G). \]
It is well known that the vertex-isoperimetric number of the
square grid $G_n$ is $n$. The down diagonal lexicographical
linear arrangement provides an example of a numbering with
vertex-isoperimetric number $n$ and the proof of
Theorem~\ref{thm:2} essentially contains a proof that this
number cannot be less than $n$. Thus in order to reduce the
bandwidth of $G_n$ we need to delete enough edges so that at
least the vertex-isoperimetric number is reduced. However, one
can reduce the vertex-isoperimetric number of $G_n$ by deleting
only one edge. For example, consider the case of $n=4$ and the
numbering given by
\[
\begin{tabular}{LLLL}
 7  & 13 & 15 & 16 \\
 4  & 8  & 12 & 14 \\
 2  & 5  &  9 & 11 \\
 1  & 3  &  6 & 10
\end{tabular}.
\]
If the edge with endpoints labelled by 7 and 13 is deleted, the
newly obtained graph has vertex-isoperimetric number 3. This
indicates that, in general, the problem of bandwidth reduction
is more difficult than the problem of vertex-isoperimetric
number reduction. It also explains why, in the course of the
proof of Theorem~\ref{thm:2}, it was relatively easy to show
that at least one edge needed to be deleted, but one had to work
harder for the second edge.

Finally, we note that the $k$-th bandwidth reduction number of
the square grid $G_n$ is bounded linearly by $2k$, for $k<n/2$,
i.e, one has to delete relatively small number of edges to
reduce the bandwidth substantially. However, it is clear that
for $k \geq n/2$ the bandwidth reduction gets more difficult. A
linear bound on the bandwidth reduction number simply cannot be
found for all $k$. After all, for $k=n$, one needs to remove all
$2n(n-1)$ edges to get to bandwidth 0.

\section*{Acknowledgments}
Thanks to Kiran Kedlaya for his input and interest.


\begin{thebibliography}{00}

\bibitem{chvatalova:bandgrid}
Jarmila Chv{\'a}talov{\'a}.
\newblock Optimal labelling of a product of two paths.
\newblock {\em Discrete Math.}, 11:249--253, 1975.

\bibitem{fishburn-w:2(n-1)}
Peter Fishburn and Paul Wright.
\newblock Bandwidth edge counts for linear arrangements of rectangular grids.
\newblock {\em J. Graph Theory}, 26(4):195--202, 1997.

\end{thebibliography}
\end{document}